\documentclass[12pt]{article}
\usepackage{mathrsfs}
\usepackage{epic,eepic,epsf,epsfig}
\usepackage{amsfonts,srcltx,mathrsfs}
\usepackage{multirow}
\usepackage{amsmath}
\usepackage{amssymb}
\usepackage{amsbsy}
\usepackage{graphicx}
\usepackage{amsfonts}
\usepackage{color}
\usepackage{setspace}

\newtheorem{lem}{Lemma}[section]%
\newtheorem{theorem}[lem]{Theorem}%

\def\nd{\mathrel{\bigm|\kern-.7em/}}

\def\f{\noindent}

\def\P\GammaL{\hbox{\rm P\GammaL}}

\def\mod{\hbox{\rm mod }}

\begin{document}
\title{Spectral conditions for graphs to contain $k$-factors}

\footnotetext{* Corresponding author}
\footnotetext{E-mails: 13021531326@163.com; zhangwq@pku.edu.cn}

\author{Xinying Tang, Wenqian Zhang*\\
{\small School of Mathematics and Statistics, Shandong University of Technology}\\
{\small Zibo, Shandong 255000, P.R. China}}
\date{}
\maketitle

\begin{abstract}
Let $G$ be a graph. The spectral radius $\rho(G)$ of $G$ is the largest eigenvalue of its adjacency matrix. For an integer $k\geq1$, a $k$-factor of $G$ is a $k$-regular spanning subgraph of $G$. Assume that $k$ and $n$ are integers satisfying $k\geq2,kn\equiv0~(\mod2)$ and $n\geq\max\left\{k^{2}+6k+7,20k+10\right\}$. Let $G$ be a graph of order $n$ and with minimum degree at least $k$. In this paper, we give a sharp lower bound of $\rho(G)$ to guarantee that $G$ contains a $k$-factor.
\bigskip

\f {\bf Keywords:} spectral radius; adjacency matrix; $k$-factor; minimum degree.

\end{abstract}

\baselineskip 17 pt

\section{Introduction}

All graphs considered in this paper are finite, undirected and simple. Let $G$ be a graph. Its vertex set and edge set are denoted by $V(G)$ and $E(G)$, respectively. Let $\overline{G}$ denote the complement of $G$. For a vertex $u$, let $d_{G}(u)$ denote its degree.  A vertex $v$ of $G$ is called a {\em neighbor} of $u$ if $uv\in E(G)$ (or $v\sim u$). Let $\delta(G)$ denote the minimum degree of $G$. For a subset $B\subseteq V(G)$, let $G[B]$ be the subgraph induced by $B$, and let $G-B$ be the graph $G[V(G)-B]$. For two vertex-disjoint subsets $W,U\subseteq V(G)$, let $e_{G}(W,U)$ be the number of edges between $W$ and $U$. For any two vertex-disjoint graphs $G_{1}$ and $G_{2}$, let $G_{1}\cup G_{2}$ be the disjoint union of them. Let $G_{1}\vee G_{2}$ be the {\em join} of $G_{1}$ and $G_{2}$, which is obtained from $G_{1}\cup G_{2}$ by connecting each vertex in $G_{1}$ to each vertex in $G_{2}$. For a positive integer $n$, let $K_{n}$  be the complete graph of order $n$.

 Let $G$ be  a graph of order $n$. Denote the vertices of $G$ by $1,2,...,n$. The {\em adjacency matrix} $A(G)$ of $G$ is an $n\times n$ matrix $(a_{ij})$, where $a_{ij}=1$ if $i\sim j$, and $a_{ij}=0$ otherwise. The {\em spectral radius} $\rho(G)$ of $G$ is the largest eigenvalue of $A(G)$.
By the Perron-Frobenius theorem, $\rho(G)$ has a non-negative eigenvector. A non-negative eigenvector corresponding to $\rho(G)$ is called a Perron vector of $G$. Moreover, if $G$ is connected, any Perron vector of $G$ is positive. For more study on this direction, one may refer to the book \cite{BH}.

Let $G$ be a graph. Assume that $0\leq a\leq b$ are integers. An $[a,b]$-factor of $G$ is a spanning subgraph with degrees between $a$ and $b$. When $a=b=k$, an $[a,b]$-factor is also called a $k$-factor. Recently, many researchers studied the spectral radius conditions for graphs to contain $[a,b]$-factors (for example, see \cite{CS,CFL,FL,FL2,FLL2,HL,HLY,O,WS}). The {\em binding number} $b(G)$ of $G$ is the minimum value of
$\frac{|N_{G}(X)|}{|X|}$ taken over all non-empty subsets $X\subseteq V(G)$ such that $N_{G}(X)\neq V(G)$, where $N_{G}(X)$ denotes the set of all the neighbors of the vertices in $X$. $G$ is called $r$-binding if $b(G)\geq r$. Very recently, Fan and Lin \cite{FL} proposed the following problem.

\medskip

\f{\bf Problem 1.} {\rm (\cite{FL})} Which 1-binding graphs $G$ with $\delta(G)\geq k$ have a $k$-factor?

\medskip

They \cite{FL} gave a spectral characterization for Problem 1 when $k=1,2$. Hao, Li and  Yu \cite{HLY} gave a spectral characterization for the bipartite analogue of Problem 1 for all $k\geq2$. In this paper, we will give a spectral characterization for Problem 1 for all $k\geq2$.

For $k\geq2$ and $n\geq3k$, let $G_{n,k}$ be the graph obtained from $K_{k}\vee(\overline{K_{k+1}}\cup K_{n-1-2k})$ by connecting one vertex in $V(\overline{K_{k+1}})$ to  $(k-1)$ vertices in $V(K_{n-1-2k})$. We can show that $G_{n,k}$ contains no $k$-factors. In fact, if $G_{n,k}$ contains a $k$-factor $F$, then there are at least $k(k+1)-(k-1)=k^{2}+1$ edges in $F$ connecting $V(\overline{K_{k+1}})$ to $V(K_{k})$. It follows that $d_{F}(u)>k$ for some vertex $u$ in $V(K_{k})$. This is obviously impossible. Our main result is the following Theorem \ref{main theo}. It is easy to see that $G_{n,k}$ is 1-binding. Thus, Theorem \ref{main theo} gives a  spectral characterization for Problem 1. Note that if a graph of order $n$ contains a $k$-factor, then $kn\equiv0~(\mod2)$.

\begin{theorem}\label{main theo}
Assume that $k\geq2,kn\equiv0~(\mod2)$ and $n\geq\max\left\{k^{2}+6k+7,20k+10\right\}$. Let $G$ be a graph of order $n$ and with $\delta(G)\geq k$. If $\rho(G)\geq\rho(G_{n,k})$, then $G$ has a $k$-factor, unless $G=G_{n,k}$.
\end{theorem}

The rest of the paper is organized as follows. In Section 2, we will include several lemmas. In Section 3, we will prove a useful lemma. In Section 4, we will give the proof of Theorem \ref{main theo}.

\section{Preliminaries}

 To prove the main results of this paper, we first include several lemmas. The following lemma is taken from \cite{BH}.

\begin{lem}{\rm (\cite{BH})}\label{subgraph}
If $H$ is a subgraph of a connected graph $G$, then $\rho(H)\leq\rho(G)$, with equality if and only if $H=G$.
\end{lem}

The following lemma is given in \cite{Zh}.

\begin{lem}{\rm (\cite{Zh})}\label{eigenvector trans}
Let $G$ be a connected graph with a Perron vector $\mathbf{x}=(x_{w})_{w\in V(G)}$. Let $u_{1}v_{1},u_{2}v_{2},\ldots,u_{s}v_{s}$ be $s\geq1$ edges of $G$, and let $a_{1}b_{1},a_{2}b_{2},\ldots,a_{t}b_{t}$ be $t\geq1$ non-edges of $G$. Let $G'$ be the graph obtained from $G$ by deleting the edges $u_{i}v_{i}$ for $1\leq i\leq s$, and adding the edges $a_{i}b_{i}$ for $1\leq i\leq t$. If $\sum_{1\leq i\leq s}x_{u_{i}}x_{v_{i}}\leq\sum_{1\leq i\leq t}x_{a_{i}}x_{b_{i}}$, and the vertex $a_{1}$ is not incident with the edges $u_{i}v_{i}$ for $1\leq i\leq s$, then $\rho(G')>\rho(G)$.
\end{lem}

The following lemma can be found in \cite{HSF}.

\begin{lem}{\rm (\cite{HSF})}\label{general bound}
Let $G$ be a graph of order $n$ and with $m$ edges satisfying $\delta(G)\geq\delta$.
Then $\rho(G)\leq \frac{\delta-1}{2}+\sqrt{2m-n\delta+\frac{(\delta+1)^{2}}{4}}$.
\end{lem}

 The following lemma is deduced from Tutte's $f$-factor theorem (see \cite{Tu}).

\begin{lem}{\rm (\cite{Tu})}\label{f-factor}
Let $G$ be a graph. For any integer $k\geq2$, $G$ has a $k$-factor if and only if for all vertex-disjoint subsets $S,T\subseteq V(G)$,
$$\delta_{G}(S,T)=\tau_{G}(S,T)+k|T|-k|S|-\sum_{u\in T}d_{G-S}(u)\leq0,$$
where $\tau_{G}(S,T)$ is the number of components $C$ of $G-(S\cup T)$ such that $e_{G}(V(C),T)+k|C|\equiv1 ~(\mod2)$. Moreover, $\delta_{G}(S,T)\equiv k|V(G)|~(\mod2)$.
\end{lem}

\section{A useful lemma}

Assume that $k\geq1$ and $n\geq3k$. Define $\mathcal{G}^{k}_{n}$ to be the set of graphs $G$ of order $n$ and with $\delta(G)\geq k$, such that there is a subset $B\subseteq V(G)$ with $|B|=k+1$ satisfying $\sum_{u\in B}d_{G}(u)\leq k^{2}+2k-1$. Clearly, $G_{n,k}\in\mathcal{G}^{k}_{n}$.

\begin{lem}\label{a lemma}
Let $\mathcal{G}^{k}_{n}$ be defined as above, where $k\geq1$ and $n\geq\frac{1}{2}k^{2}+3k+1$. Then $G_{n,k}$ is the unique extremal graph with the maximum spectral radius in $\mathcal{G}^{k}_{n}$.
\end{lem}

\f{\bf Proof:} Let $G$ be an extremal graph with the maximum spectral radius in $\mathcal{G}^{k}_{n}$. It suffices to prove that $G=G_{n,k}$. Let $B$ be a subset of $V(G)$ with $|B|=k+1$, such that $\sum_{u\in B}d_{G}(u)\leq k^{2}+2k-1$. Clearly, $G-B$ is a complete graph and $\sum_{u\in B}d_{G}(u)= k^{2}+2k-1$ by Lemma \ref{subgraph}, since $G$ has the maximum spectral radius. Let $\rho=\rho(G)$.
By Lemma \ref{subgraph} again, we have $\rho>\rho(K_{n-1-k})=n-2-k$, since $G$ contains $K_{n-1-k}$ as a proper subgraph. Let $\mathbf{x}=(x_{u})_{u\in V(G)}$ be a Perron vector of $G$.

\medskip

\f{\bf Claim 1.} For any graph $G'\in\mathcal{G}^{k}_{n}$, $\rho(G')<n-1-k$.

\medskip

\f{\bf Proof of Claim 1.}  Let $B'$ be a subset of $V(G')$ with $|B'|=k+1$, such that
$$\sum_{u\in B'}d_{G'}(u)\leq k^{2}+2k-1.$$
Then
$$e(G')\leq(\sum_{u\in B'}d_{G'}(u))+e(G'-B')\leq k^{2}+2k-1+\frac{(n-1-k)(n-2-k)}{2}.$$
Recall that $\delta(G')\geq k$. By Lemma \ref{general bound}, we have
 \begin{equation}
\begin{aligned}
\rho(G')
&\leq\frac{k-1}{2}+\sqrt{2e(G')-nk+\frac{(k+1)^{2}}{4}}\\
&\leq\frac{k-1}{2}+\sqrt{2k^{2}+4k-2+(n-1-k)(n-2-k)-nk+\frac{(k+1)^{2}}{4}}\\
&=\frac{k-1}{2}+\sqrt{n^{2}-(3k+3)n+\frac{21}{4}k^{2}+\frac{15}{2}k+\frac{1}{4}}\\
&<n-1-k~(since~n\geq\frac{1}{2}k^{2}+3k+1).
\end{aligned}\notag
\end{equation}
This finishes the proof Claim 1. \hfill$\Box$

\medskip

\f{\bf Claim 2.} $G$ can be obtained from $K_{k}\vee(\overline{K_{k+1}}\cup K_{n-1-2k})$ by adding $(k-1)$ edges between $V(\overline{K_{k+1}})$ and $V(K_{n-1-2k})$.

\medskip

\f{\bf Proof of Claim 2.} Since $\delta(G)\geq k$ and $\sum_{u\in B}d_{G}(u)= k^{2}+2k-1$, we have $d_{G}(u)\leq 2k-1$ for any $u\in B$.
Denote $B=\left\{u_{1},u_{2},...,u_{k+1}\right\}$ and $V(G)-B=\left\{v_{1},v_{2},...,v_{n-1-k}\right\}$.  Without loss of generality, assume that $x_{u_{1}}\geq x_{u_{2}}\geq\cdots\geq x_{u_{k+1}}$ and $x_{v_{1}}\geq x_{v_{2}}\geq\cdots\geq x_{v_{n-1-k}}$.
Now we prove that $x_{v_{n-1-k}}\geq x_{u_{1}}$. Let $d_{G}(u_{1})=d$, and let $r$ be the number of neighbors of $u_{1}$ in $B$.
Since
$$\rho x_{u_{1}}=(\sum_{u\in B,u\sim u_{1}}x_{u})+(\sum_{u\in V(G)-B,u\sim u_{1}}x_{u})\leq rx_{u_{1}}+\sum_{1\leq i\leq d-r}x_{v_{i}},$$
we obtain
$$x_{u_{1}}\leq\frac{1}{\rho-r}\sum_{1\leq i\leq d-r}x_{v_{i}}.$$
Since
\begin{equation}
\begin{aligned}
\rho x_{v_{n-1-k}}&=(\sum_{u\in B,u\sim v_{n-1-k}}x_{u})+(\sum_{u\in V(G)-B,u\sim v_{n-1-k}}x_{u})\\
&\geq\sum_{u\in V(G)-B,u\sim v_{n-1-k}}x_{u}\\
&\geq(n-2-k-d+r)x_{v_{n-1-k}}+\sum_{1\leq i\leq d-r}x_{v_{i}},
\end{aligned}\notag
\end{equation}
we obtain
$$x_{v_{n-1-k}}\geq\frac{1}{\rho-(n-2-k)+d-r}\sum_{1\leq i\leq d-r}x_{v_{i}}.$$
Note that $\rho-(n-2-k)+d-r\leq\rho-r$, since $n\geq\frac{1}{2}k^{2}+3k+1\geq3k+1\geq k+d+2$. It follows that $x_{v_{n-1-k}}\geq x_{u_{1}}$.

Now we show that there are no edges inside $B$. Suppose not. Without loss of generality, assume that $u_{1}u_{2}$ is an edge in $G$. Since $n\geq\frac{1}{2}k^{2}+3k+1$, there is a vertex $v\in V(G)-B$ such that $v$ is not adjacent to $u_{1}$ and $u_{2}$. Let $G_{1}$ be the graph obtained from $G$ by deleting the edge $u_{1}u_{2}$ and adding the edges $vu_{1}$ and $vu_{2}$. Clearly, $G_{1}$ is in $\mathcal{G}^{k}_{n}$. By Lemma \ref{eigenvector trans}, noting that $x_{v}\geq x_{u_{1}}$, we have $\rho(G_{1})>\rho(G)$, which contradicts the choice of $G$. Hence, there are no edges inside $B$. Since $\delta(G)\geq k$, using a similar discussion, we can show that $u_{i}$ is adjacent to $v_{1},v_{2},...,v_{k}$ for any $1\leq i\leq k+1$. Hence, $G$ can be obtained from $K_{k}\vee(\overline{K_{k+1}}\cup K_{n-1-2k})$ by adding $(k-1)$ edges between $V(\overline{K_{k+1}})$ and $V(K_{n-1-2k})$. This finishes the proof Claim 2. \hfill$\Box$

Denote $V(K_{k})=\left\{w_{1},w_{2},...,w_{k}\right\},V(\overline{K_{k+1}})=\left\{u_{1},u_{2},...,u_{k+1}\right\}$ and $V(K_{n-1-2k})=\left\{v_{1},v_{2},...,v_{n-1-2k}\right\}$. By Claim 2, $G$ is obtained from $K_{k}\vee(\overline{K_{k+1}}\cup K_{n-1-2k})$ by adding $(k-1)$ edges between $\left\{u_{1},u_{2},...,u_{k+1}\right\}$ and $\left\{v_{1},v_{2},...,v_{n-1-2k}\right\}$. Without loss of generality, assume that $x_{u_{1}}\geq x_{u_{2}}\geq\cdots\geq x_{u_{k+1}}$ and $x_{v_{1}}\geq x_{v_{2}}\geq\cdots\geq x_{v_{n-1-2k}}$.
By symmetry, we see $x_{w_{1}}= x_{w_{2}}=\cdots= x_{w_{k}}$.

Let $s\geq k$ be the largest integer such that $v_{s}$ is adjacent to $u_{1}$ in $G$. We can show that $v_{i}$ is adjacent to $u_{1}$ for any $1\leq i\leq s$ in $G$. Otherwise assume that $v_{j}u_{1}$ is not an edge of $G$ for some $1\leq j< s$. Let $G_{2}$ be the graph obtained from $G$ by deleting the edge $v_{s}u_{1}$ and adding the edge $v_{j}u_{1}$. Clearly, $G_{2}$ is in $\mathcal{G}^{k}_{n}$. By Lemma \ref{eigenvector trans}, noting that $x_{v_{j}}\geq x_{v_{s}}$, we have $\rho(G_{2})>\rho(G)$, which contradicts the choice of $G$. Thus, $v_{i}u_{1}$ is an edge for any  $1\leq i\leq s$.  We can show that $v_{i}u_{\ell}$ is not an edge of $G$ for any $i>s$ and $1\leq \ell\leq k+1$. In fact, if $v_{i_{0}}u_{\ell_{0}}$ is an edge of $G$ for some $i_{0}>s$ and $1\leq \ell_{0}\leq k+1$, then $\ell_{0}\geq2$ by the choice of $s$. Let $G_{3}$ be the graph obtained from $G$ by deleting the edge $v_{i_{0}}u_{\ell_{0}}$ and adding the edge $v_{i_{0}}u_{1}$. Clearly, $G_{3}$ is in $\mathcal{G}^{k}_{n}$. By Lemma \ref{eigenvector trans}, noting that $x_{u_{\ell_{0}}}\leq x_{u_{1}}$, we have $\rho(G_{3})>\rho(G)$, which contradicts the choice of $G$. Thus we obtain that $v_{i}u_{\ell}$ is not an edge for any $i>s$ and $1\leq \ell\leq k+1$.

 By symmetry, we have $x_{v_{s+1}}=x_{v_{i}}$ for any $s+2\leq i\leq n-1-2k$. By $A(G)\mathbf{x}=\rho\mathbf{x}$, we have
 $$\rho x_{v_{s+1}}=kx_{w_{1}}+x_{v_{1}}+(\sum_{2\leq i\leq s}x_{v_{i}})+(n-2-2k-s)x_{v_{s+1}}\geq (k+1)x_{v_{1}}+(n-3-2k)x_{v_{s+1}}.$$
It follows that
$$x_{v_{s+1}}\geq\frac{(k+1)x_{v_{1}}}{\rho+3+2k-n}.$$
Since $\rho<n-1-k$ by Claim 1, we see $$\frac{x_{v_{s+1}}}{x_{v_{1}}}\geq\frac{k+1}{\rho+3+2k-n}>\frac{k+1}{k+2}.$$

Recall that $G_{n,k}$ is the graph obtained from $K_{k}\vee(\overline{K_{k+1}}\cup K_{n-1-2k})$ by connecting the vertex $u_{1}$ to the vertices $v_{1},v_{2},...,v_{k-1}$. Now we prove that $G=G_{n,k}$. If $s=k-1$, then $G=G_{n,k}$, as desired. Now assume that $1\leq s<k-1$, implying $k>2$. We will obtain a  contradiction. Denote $\rho(G_{n,k})=\rho'$. Then $\rho'<n-1-k$ by Claim 1. Let $\mathbf{y}=(y_{u})_{u\in V(G_{n,k})}$ be a Perron vector of $G_{n,k}$. By symmetry, we have $y_{u_{2}}=y_{u_{3}}=\cdots=y_{u_{k+1}}$, $y_{w_{1}}=y_{w_{2}}=\cdots=y_{w_{k}}$ and $y_{v_{1}}=y_{v_{2}}=\cdots=y_{v_{k-1}}\geq y_{v_{k}}=y_{v_{k+1}}=\cdots=y_{v_{n-1-2k}}$.
By $A(G_{n,k})\mathbf{y}=\rho'\mathbf{y}$, we have
$$\rho' y_{u_{1}}=ky_{w_{1}}+(k-1)y_{v_{1}},$$
and
$$\rho' y_{u_{2}}=ky_{w_{1}}.$$
 Since $$\rho' y_{v_{k}}=ky_{w_{1}}+(k-1)y_{v_{1}}+(n-1-3k)y_{v_{k}}\geq ky_{w_{1}}+(n-2-2k)y_{v_{k}},$$
 we obtain that $$ y_{v_{1}}\geq y_{v_{k}}\geq\frac{ky_{w_{1}}}{\rho'+2+2k-n}.$$
Then $$\rho' y_{u_{1}}=ky_{w_{1}}+(k-1)y_{v_{1}}\geq ky_{w_{1}}+(k-1)\frac{ky_{w_{1}}}{\rho'+2+2k-n}.$$
Hence, $$\frac{y_{u_{1}}}{y_{u_{2}}}=\frac{\rho'y_{u_{1}}}{\rho'y_{u_{2}}}\geq1+\frac{k-1}{\rho'+2+2k-n}=\frac{\rho'+1+3k-n}{\rho'+2+2k-n}.$$
Let $\mathbf{x}^{T}$ denote the transpose of $\mathbf{x}$.  Thus,
\begin{equation}
\begin{aligned}
&(\rho'-\rho)\mathbf{x}^{T}\mathbf{y}\\
&=\mathbf{x}^{T}(A(G_{n,k})-A(G))\mathbf{y}\\
&=(\sum_{u_{1}v_{i}\in (E(G_{n,k})-E(G))}(x_{u_{1}}y_{v_{i}}+x_{v_{i}}y_{u_{1}}))-(\sum_{u_{i}v_{j}\in (E(G)-E(G_{n,k}))}(x_{u_{i}}y_{v_{j}}+x_{v_{j}}y_{u_{i}}))\\
&\geq(k-1-s)(x_{u_{1}}y_{v_{1}}+x_{v_{s+1}}y_{u_{1}}-x_{u_{2}}y_{v_{1}}-x_{v_{1}}y_{u_{2}})\\
&\geq(k-1-s)(x_{v_{s+1}}y_{u_{1}}-x_{v_{1}}y_{u_{2}})\\
&>(k-1-s)x_{v_{1}}y_{u_{2}}(\frac{k+1}{k+2}\frac{\rho'+1+3k-n}{\rho'+2+2k-n}-1)\\
&=(k-1-s)x_{v_{1}}y_{u_{2}}\frac{n-\rho'+k^{2}-2k-3}{(k+2)(\rho'+2+2k-n)}\\
&>(k-1-s)x_{v_{1}}y_{u_{2}}\frac{k^{2}-k-2}{(k+2)(\rho'+2+2k-n)}~(since ~\rho'<n-1-k)\\
&\geq0.
\end{aligned}\notag
\end{equation}
That is $(\rho'-\rho)\mathbf{x}^{T}\mathbf{y}>0$, implying that $\rho'>\rho$. But this contradicts the choice of $G$.
This completes the proof. \hfill$\Box$

\section{Proof of Theorem \ref{main theo}}

Let $G$ be a graph. For two vertex-disjoint subsets $S,T\subseteq V(G)$, let $q_{G}(S,T)$ denote the number of components of $G-(S\cup T)$. For positive integers $n,k,$ let $\mathcal{G}_{n,k}$ be the set of graphs $G$ of order $n$ with $\delta(G)\geq k$, such that there are two vertex-disjoint subsets $S,T\subseteq V(G)$ satisfying
$$\sum_{u\in T}d_{G-S}(u)\leq k|T|-k|S|-2+q_{G}(S,T).$$

\begin{lem}\label{main lemma}
For $k\geq2$ and $n\geq\max\left\{k^{2}+6k+7,20k+10\right\}$, $G_{n,k}$ is the unique extremal graph with the maximum spectral radius in  $\mathcal{G}_{n,k}$.
\end{lem}

\f{\bf Proof:} Recall that $G_{n,k}$ is obtained from $K_{k}\vee(\overline{K_{k+1}}\cup K_{n-1-2k})$ by adding $k-1$ edges connecting one vertex of $\overline{K_{k+1}}$ to $(k-1)$ vertices of $K_{n-1-2k}$. It is easy to check that $G_{n,k}\in \mathcal{G}_{n,k}$ by letting $S=V(K_{k})$ and $T=V(\overline{K_{k+1}})$. Let $G$ be an extremal graph with the maximum spectral radius in  $\mathcal{G}_{n,k}$. It suffices to prove that $G=G_{n,k}$.

Since $G_{n,k}$ contains $K_{n-1-k}$ as a proper subgraph, we have $\rho(G_{n,k})>\rho(K_{n-1-k})=n-2-k$ by Lemma \ref{subgraph}. Then $\rho(G)\geq\rho(G_{n,k})>n-2-k$. By Lemma \ref{general bound}, we have
 $$\rho(G)\leq\frac{k-1}{2}+\sqrt{2e(G)-kn+\frac{(k+1)^{2}}{4}}.$$
Noting that $\rho(G)>n-k-2$, we obtain that
$$e(G)>\frac{1}{2}n^{2}-(k+\frac{3}{2})n+(k+1)^{2},$$
and thus
\begin{equation}
\begin{aligned}
e(\overline{G})<(k+1)n-(k+1)^{2}.
\end{aligned}
\end{equation}

Since $G\in\mathcal{G}_{n,k}$, there are two vertex-disjoint subsets $S,T\subseteq V(G)$ satisfying
$$\sum_{u\in T}d_{G-S}(u)\leq k|T|-k|S|-2+q_{G}(S,T).$$
We can choose such $S$ and $T$ that $|S\cup T|$ is maximum. Set $s=|S|,t=|T|$ and $q=q_{G}(S,T)$.
Then \begin{equation}
\begin{aligned}
\sum_{u\in T}d_{G-S}(u)\leq kt-ks-2+q.
\end{aligned}
\end{equation}
Let  $Q_{1},Q_{2},...,Q_{q}$ be the  components of $G-(S\cup T)$, where $n_{i}=|Q_{i}|$ for $1\leq i\leq q$. Without loss of generality, assume that $n_{1}\geq n_{2}\geq\cdots\geq n_{q}$.
Let $\mathbf{x}=(x_{u})_{u\in V(G)}$ be a Perron vector of $G$.

\medskip

\f{\bf Claim 1.} $d_{G}(u)=n-1$ for any $u\in S$, and $Q_{i}$ is a complete graph for each $1\leq i\leq q$.

\medskip

\f{\bf Proof of Claim 1.} Let $G_{1}$ be the graph obtained from $G$ by adding edges such that $d_{G_{1}}(u)=n-1$ for any $u\in S$, and $Q_{i}$ is a complete graph for each $1\leq i\leq q$. By Lemma \ref{subgraph}, $\rho(G)\leq\rho(G_{1})$ with equality if and only if $G=G_{1}$. Clearly,  $\delta(G_{1})\geq k$. Note that $q_{G_{1}}(S,T)=q_{G}(S,T)$ and $\sum_{u\in T}d_{G_{1}-S}(u)=\sum_{u\in T}d_{G-S}(u)$. It follows that
 $$\sum_{u\in T}d_{G_{1}-S}(u)\leq k|T|-k|S|-2+q_{G_{1}}(S,T).$$
Then $G_{1}\in\mathcal{G}_{n,k}$. Since $G$ is a spanning subgraph of $G_{1}$, we see $G=G_{1}$ by the choice of $G$.  This finishes the proof of Claim 1. \hfill$\Box$

\medskip

\f{\bf Claim 2.} For any $i\geq1$, if $Q_{i}$ is a singleton $\left\{w_{i}\right\}$, then $d_{G-S}(w_{i})=k$. If $n_{i}\geq2$, then $d_{G-S}(v)\geq k+1$ for any $v\in V(Q_{i})$.

\medskip

\f{\bf Proof of Claim 2.} First assume that $Q_{i}$ is a singleton $\left\{w_{i}\right\}$. We will show $d_{G-S}(w_{i})=k$. If $d_{G-S}(w_{i})\leq k-1$, let $S'=S$ and $T'=T\cup\left\{w_{i}\right\}$.  Clearly, $q_{G}(S',T')= q_{G}(S,T)-1, \sum_{u\in T'}d_{G-S'}(u)\leq k-1+\sum_{u\in T}d_{G-S}(u)$ and  $|S'\cup T'|=|S\cup T|+1$. It follows  that
$$\sum_{u\in T'}d_{G-S'}(u)\leq k|T'|-k|S'|-2+q_{G}(S',T').$$
This contradicts the choices of $S$ and $T$, since $|S\cup T|$ is maximum. If $d_{G-S}(w_{i})\geq k+1$, let $S'=S\cup\left\{w_{i}\right\}$ and $T'=T$. Clearly, $q_{G}(S',T')= q_{G}(S,T)-1, \sum_{u\in T'}d_{G-S'}(u)\leq -(k+1)+\sum_{u\in T}d_{G-S}(u)$ and  $|S'\cup T'|=|S\cup T|+1$. It follows  that $$\sum_{u\in T'}d_{G-S'}(u)\leq k|T'|-k|S'|-2+q_{G}(S',T').$$
This is a contradiction, since $|S\cup T|$ is maximum. Consequently, $d_{G-S}(w_{i})=k$.

Now assume that $n_{i}\geq2$. Let $v\in V(Q_{i})$. We will show $d_{G-S}(v)\geq k+1$. If $d_{G-S}(v)\leq k$, let $S'=S$ and $T'=T\cup\left\{v\right\}$.  Clearly, $q_{G}(S',T')= q_{G}(S,T), \sum_{u\in T'}d_{G-S'}(u)\leq k+\sum_{u\in T}d_{G-S}(u)$ and  $|S'\cup T'|=|S\cup T|+1$. It follows  that $$\sum_{u\in T'}d_{G-S'}(u)\leq k|T'|-k|S'|-2+q_{G}(S',T').$$
This is a contradiction, since $|S\cup T|$ is maximum. Hence, $d_{G-S}(v)\geq k+1$.
 This finishes the proof of Claim 2. \hfill$\Box$

\medskip

\f{\bf Claim 3.} For any $i\geq2$, we have $n_{i}\leq k-1$.

\medskip

\f{\bf Proof of Claim 3.} Suppose $n_{i}\geq k$.  Then $n_{1}\geq n_{i}\geq k$. By Claim 2, $d_{G-S}(u)\geq k+1$ for any $u\in V(Q_{1})\cup V(Q_{i})$. Without loss of generality, assume that $\sum_{u\in V(Q_{1})}x_{u}\geq\sum_{u\in V(Q_{i})}x_{u}$. Let $v$ be a vertex in $V(Q_{i})$.
Let $G_{2}$ be the graph obtained from $G$ by deleting the edges between $v$ and $V(Q_{i})-\left\{v\right\}$, and adding the edges between $v$ and $V(Q_{1})$. Clearly, $\delta(G_{2})\geq k$, $q_{G_{2}}(S,T)= q_{G}(S,T)$ and $\sum_{u\in T}d_{G_{2}-S}(u)= \sum_{u\in T}d_{G-S}(u)$. It follows  that $$\sum_{u\in T}d_{G_{2}-S}(u)\leq k|T|-k|S|-2+q_{G_{2}}(S,T).$$
Hence, $G_{2}\in \mathcal{G}_{n,k}$. But $\rho(G_{2})>\rho(G)$ by Lemma \ref{eigenvector trans}, which contradicts the choice of $G$.
Thus, $n_{i}\leq k-1$.  This finishes the proof of Claim 3. \hfill$\Box$

\medskip

\f{\bf Claim 4.} $e_{G}(T,V(Q_{i}))\geq1$ for any $2\leq i\leq q$. Consequently, $s\leq t-1$.

\medskip

\f{\bf Proof of Claim 4.} Assume that $2\leq i\leq q$. If $n_{i}=1$, there are $k$ edges between $T$ and $V(Q_{i})$ by Claim 2.
 If $n_{i}\geq2$, then $n_{i}\leq k-1$ by Claim 3, and $d_{G-S}(u)\geq k+1$ for any $u\in V(Q_{i})$ by Claim 2. This implies that
  there is at least one edge between $u$ and $T$. Hence, $e_{G}(T,V(Q_{i}))\geq1$ for any $2\leq i\leq q$.
  It follows that $\sum_{u\in T}d_{G-S}(u)\geq q-1$.
 Recall that
$$\sum_{u\in T}d_{G-S}(u)\leq kt-ks-2+q$$ in $(2)$.
Then $q-1\leq kt-ks-2+q$, implying that $s\leq t-1$. This finishes the proof of Claim 4. \hfill$\Box$

\medskip

\f{\bf Claim 5.} $t\leq\frac{1}{2}n-3k$ and $q\geq1$.

\medskip

\f{\bf Proof of Claim 5.} We first show $t\leq\frac{1}{2}n-3k$. Suppose that $t>\frac{1}{2}n-3k$. Considering the non-edges inside $T$ and among the components $Q_{1},Q_{2},...,Q_{q}$, we have
\begin{equation}
\begin{aligned}
e(\overline{G})
&\geq\frac{t(t-1)}{2}-\frac{1}{2}(\sum_{u\in T}d_{G-S}(u))+q-1\\
&\geq\frac{t(t-1)}{2}-\frac{1}{2}kt~(using~(2))\\
&\geq\frac{(\frac{1}{2}n-3k)(\frac{1}{2}n-3k-1)}{2}-\frac{1}{2}k(\frac{1}{2}n-3k)\\
&=\frac{1}{8}n^{2}-\frac{1}{4}(7k+1)n+6k^{2}+\frac{3}{2}k\\
&\geq(k+1)n-(k+1)^{2}~(since~n\geq20k+10),
\end{aligned}\notag
\end{equation}
which contradicts the formula $(1)$. Hence, $t\leq\frac{1}{2}n-3k$.

If $q=0$, then $n=s+t$. Since $s\leq t-1$ by Claim 4, we have $t\geq\frac{1}{2}n$. But this contradicts the proved result $t\leq\frac{1}{2}n-3k$.
 Hence, $q\geq1$. This finishes the proof of Claim 5. \hfill$\Box$

\medskip

Depending on the value of $q$, we have the following 3 cases to handle  by Claim 5.

\medskip

\f{\bf Case 1.} $q=1$.

\medskip

Now $(2)$ becomes
$$\sum_{u\in T}d_{G-S}(u)\leq kt-ks-1.$$
Since $\sum_{u\in T}d_{G-S}(u)\geq(k-s)t$ as $\delta(G)\geq k$, we see $(k-s)t\leq kt-ks-1$, and thus $(k-t)s\leq-1$. This implies that $t\geq k+1$.
Considering the non-edges between $T$ and $V(G)-(S\cup T)$, we have
\begin{equation}
\begin{aligned}
e(\overline{G})
&\geq|T||V(G)-(T\cup S)|-\sum_{u\in T}d_{G-S}(u)\\
&\geq t(n-s-t)-k(t-s)+1.
\end{aligned}
\end{equation}

\f{\bf Subcase 1.1.} $t=k+1$.

Note that $s\leq t-1=k$. Then $$\sum_{u\in T}d_{G}(u)\leq ts+\sum_{u\in T}d_{G-S}(u)\leq k^{2}+2k-1.$$ Hence, $G\in \mathcal{G}^{k}_{n}$. By Lemma \ref{a lemma}, we have $\rho(G)\leq\rho(G_{n,k})$ with equality if and only if $G=G_{n,k}$. Hence, $G=G_{n,k}$ by the choice of $G$.

\medskip

\f{\bf Subcase 1.2.} $k+2\leq t\leq\frac{n}{2}-3k$.

\medskip

Recall that $s\leq t$ by Claim 4.
By $(3)$ we see
\begin{equation}
\begin{aligned}
e(\overline{G})&\geq t(n-s-t)-k(t-s)+1\\
&=-(t-k)s+t(n-t)-kt+1\\
&\geq-(t-k)t+t(n-t)t-kt+1~(since~t\geq k+2)\\
&=t(n-2t)+1\\
&\geq(k+2)(n-2(k+2))+1~(since~n-2t\geq2(k+2))\\
&\geq(k+1)n-(k+1)^{2}~(since~n\geq k^{2}+6k+7),
\end{aligned}\notag
\end{equation}
a contradiction by $(1)$.

\medskip

\f{\bf Case 2.} $q=2$.

\medskip

Now $(2)$ becomes
$$\sum_{u\in T}d_{G-S}(u)\leq kt-ks.$$
Let $C=T\cup V(Q_{2})$. If $n_{2}=1$, then $t\geq k$ by Claim 2. If $n_{2}\geq2$, then $n_{2}\leq k-1$ by Claim 3, and $|C|\geq k+2$ by Claim 2. Hence, $|C|\geq k+1$ with equality only if $n_{2}=1$ and $t=k$. In either case, $|C|\leq\frac{1}{2}n-2k-1$ since $t\leq\frac{1}{2}n-3k$ by Claim 5.
Considering the non-edges between $C$ and $V(G)-(S\cup C)$, we have
\begin{equation}
\begin{aligned}
e(\overline{G})
&\geq|C||V(G)-(S\cup C)|-\sum_{u\in T}d_{G-S}(u)\\
&\geq|C|(n-s-|C|)-k(t-s).
\end{aligned}
\end{equation}

\medskip

\f{\bf Subcase 2.1.} $|C|=k+1$.

\medskip

Then $n_{2}=1$ and $t=k$. Thus, $s\leq k-1$ by Claim 4.
Then
 $$\sum_{u\in C}d_{G}(u)=d_{G}(w_{2})+\sum_{u\in T}d_{G}(u)\leq (k+s)+(ts+\sum_{u\in T}d_{G-S}(u))\leq k^{2}+2k-1.$$ Hence, $G\in \mathcal{G}^{k}_{n}$. By Lemma \ref{a lemma}, we have $\rho(G)\leq\rho(G_{n,k})$ with equality if and only if $G=G_{n,k}$. Hence, $G=G_{n,k}$ by the choice of $G$.

\medskip

\f{\bf Subcase 2.2.} $k+2\leq |C|\leq\frac{n}{2}-2k-1$.

\medskip

Recall that $s\leq t$ by Claim 4.
By $(4)$ we have
\begin{equation}
\begin{aligned}
e(\overline{G})&\geq|C|(n-s-|C|)-k(t-s)\\
&=-(|C|-k)s+|C|(n-|C|)-kt\\
&\geq-(|C|-k)t+|C|(n-|C|)-kt\\
&=|C|(n-t-|C|)\\
&\geq|C|(n-2|C|)\\
&\geq(k+2)(n-2(k+2))~(since~n-2|C|\geq2(k+2))\\
&\geq(k+1)n-(k+1)^{2}~(since~n\geq k^{2}+6k+7),
\end{aligned}\notag
\end{equation}
a contradiction by $(1)$.

\medskip

\f{\bf Case 3.} $q\geq3$.

\medskip

By Claim 4, there is at least one edge between $T$ and $V(Q_{i})$ for any $2\leq i\leq q$. If $q=3$, define $H=G$. If $q\geq4$, define $H$ to be the graph obtained from $G$ by deleting one edge between $V(Q_{i})$ and $T$ for any $i\geq4$, and connecting $V(Q_{i})$ to $V(Q_{1})$ for any $i\geq4$. Clearly, $e(H)\geq e(G)$ in either case, implying
$$e(\overline{H})\leq e(\overline{G})< (k+1)n-(k+1)^{2}.$$
 Moreover, $H-(S\cup T)$ has 3 components $Q'_{1},Q'_{2},Q'_{3}$ satisfying $V(Q'_{1})=V(Q_{1})\cup (\cup_{4\leq i\leq q}V(Q_{i}))$, $Q'_{2}=Q_{2}$ and $Q'_{3}=Q_{3}$. That is $q_{H}(S,T)=3$.
 Note that
$$\sum_{u\in T}d_{H-S}(u)= -(q-3)+\sum_{u\in T}d_{G-S}(u).$$
Thus,
$$\sum_{u\in T}d_{H-S}(u)\leq kt-ks-2+q_{H}(S,T)=k(t-s)+1.$$

 Let $D=T\cup V(Q_{2})\cup V(Q_{3})$. It is easy to see that $|D|\geq k+2$ by Claim 2. Since $n_{2},n_{3}\leq k-1$ by Claim 3 and $t\leq\frac{1}{2}n-3k$ by Claim 5, we have $|D|\leq\frac{1}{2}n-k-2$.
Thus, $$k+2\leq|D|\leq\frac{1}{2}n-k-2.$$
Recall $s\leq t$ by Claim 4. Considering the non-edges of $H$ between $D$ and $V(G)-(S\cup D)$, and one non-edge between $Q_{2}$ and $Q_{3}$, we have
\begin{equation}
\begin{aligned}
e(\overline{H})
&\geq|D||V(G)-(S\cup D)|+1-\sum_{u\in T}d_{H-S}(u)\\
&\geq|D|(n-s-|D|)-k(t-s)\\
&=-(|D|-k)s+|D|(n-|D|)-kt\\
&\geq-(|D|-k)t+|D|(n-|D|)-kt\\
&=|D|(n-t-|D|)\\
&\geq |D|(n-2|D|)\\
&\geq(k+2)(n-2(k+2)) ~(since~n-2|D|\geq2(k+2))\\
&\geq(k+1)n-(k+1)^{2}~(since~n\geq k^{2}+6k+7),
\end{aligned}\notag
\end{equation}
a contradiction by $(1)$.
This completes the proof. \hfill$\Box$

\medskip

\f{\bf The proof of Theorem \ref{main theo}.} Let $G$ be a graph of order $n$ and with $\delta(G)\geq k$, such that $\rho(G)\geq\rho(G_{n,k})$ and $G$ contains no $k$-factors. It suffices to prove that $G=G_{n,k}$.
Since $G$ has no $k$-factors, by Lemma \ref{f-factor}, there are two  vertex-disjoint subsets $S,T\subseteq V(G)$, such that
$$\delta_{G}(S,T)=\tau_{G}(S,T)+k|T|-k|S|-\sum_{u\in T}d_{G-S}(u)>0,$$
where $\tau_{G}(S,T)$ is the number of components $C$ of $G-(S\cup T)$ such that $e_{G}(V(C),T)+k|C|\equiv1 ~(\mod2)$. Moreover, $\delta_{G}(S,T)\equiv kn~(\mod2)$. Since $kn$ is even by assumption, we have $\delta_{G}(S,T)\geq2$. Then
$$\tau_{G}(S,T)+k|T|-k|S|-\sum_{u\in T}d_{G-S}(u)\geq2.$$
Recall that $q_{G}(S,T)$ is the number of components of $G-(S\cup T)$. Clearly, $q_{G}(s,t)\geq\tau_{G}(S,T)$.
Thus,  $$\sum_{u\in T}d_{G-S}(u)\leq k|T|-k|S|-2+q_{G}(S,T).$$
This implies that $G\in \mathcal{G}_{n,k}$. By Lemma \ref{main lemma}, we have $\rho(G)\leq\rho(G_{n,k})$ with equality if and only if $G=G_{n,k}$.
Hence, $G=G_{n,k}$ by the choice of $G$. This completes the proof. \hfill$\Box$

\medskip

\f{\bf Data availability statement}

\medskip

There is no associated data.

\medskip

\f{\bf Declaration of Interest Statement}

\medskip

There is no conflict of interest.

\end{document}